\documentclass[12pt, twoside, leqno]{article}

\usepackage{amsmath,amsthm}
\usepackage{amssymb}

\usepackage{enumerate}

\usepackage{graphicx}

\newtheorem{thm}{Theorem}[section]
\newtheorem{cor}[thm]{Corollary}
\newtheorem{lem}[thm]{Lemma}
\newtheorem{prop}[thm]{Proposition}

\theoremstyle{definition}

\numberwithin{equation}{section}

\usepackage{pgf,tikz}
\usetikzlibrary{arrows}

\begin{document}

\baselineskip=13pt


\title{Some algebraic properties of bipartite Kneser graphs}

\author{S. Morteza Mirafzal,  Ali Zafari\\
Department of Mathematics \\
  Lorestan University, Khorramabad, Iran\\
E-mail: mirafzal.m@lu.ac.ir\\
 E-mail: zafari.math.pu@gmail.com}

\date{}

\maketitle


\renewcommand{\thefootnote}{}



\renewcommand{\thefootnote}{\arabic{footnote}}
\setcounter{footnote}{0}


\begin{abstract}
Let $n$ and $k$ be   integers with  $n> k\geq1$ and $[n] = \{1, 2, ... , n\} $.  The $bipartite \ Kneser \ graph$ $H(n, k)$ is the
graph with the all $k$-element and  all ($n-k$)-element subsets of $[n]  $ as vertices,
and there is an edge between any two vertices, when one is a subset of the other.
In this paper, we show   that $H(n, k)$  is an arc-transitive graph. Also, we show that $H(n,1)$ is a distance-transitive  Cayley graph. Finally,  we determine the automorphism group  of the graph $H(n, 1)$ and show that $Aut(H(n, 1)) \cong Sym([n] ) \newline  \times \mathbb{Z}_2$,  where $\mathbb{Z}_2$ is the cyclic group of order $2$. Moreover,  we pose some open problems about the automorphism group of the  bipartite Kneser graph $H(n, k)$.     \
\end{abstract}
\section{Introduction}
For a positive integer $n >  1 $,  let $[n]  = \{1, 2, ... , n\} $   and $V$ be the set of all $k$-subsets
and $(n-k)$-subsets of $[n]$.
 The $bipartite\ Kneser \newline  graph$  $H(n, k)$ has
$V$ as its vertex set, and two vertices $A,  B$ are connected if and only if $A \subset B$ or $B\subset A$. If $n = 2k$ it is obvious that
we do not have any edges, and in such a case, $H(n, k)$ is a null graph, and  hence we assume that  $n \geq 2k + 1$.
It follows from the definition of the graph $H(n, k)$ that it has
 2${n}\choose{k}$  vertices and the degree of each of its  vertices   is
 ${n-k}\choose{k}$= ${n-k}\choose{n-2k}$, hence it is a regular graph. It is clear that $H(n, k)$ is a bipartite graph.
In fact,
if  $V_1=\{ v\in V(H(n ,k))   |  \  |v| =k \}$ and $V_2=\{ v\in V(H(n ,k))  | \  |v| =n-k \}$, then $\{ V_1, V_2\}$
is a partition of $V(H(n ,k))$ and every edge of $H(n, k)$ has a vertex in $V_1$ and a  vertex in $V_2$ and
$| V_1 |=| V_2 |$.  It is an easy task to show that   the graph $H(n, k)$  is a connected graph. The bipartite Kneser graph $H(2n-1, n-1)$ is known as the middle cube $MQ_n$ (Dalfo,   Fiol,   Mitjana [2]) or regular hyper-star graph $HS(2n,n)$(Kim,   Cheng,   Liptak,   Lee [5]). It was conjectured by Dejter, Erdos, and Havel [6] among others, that $MQ_n$
is Hamiltonian.  The graph $MQ_n$ has been studied by various authors  [2,5,6,8,12].
 Recently, Mutze [12]  showed that the bipartite Kneser graph $H(n, k)$ has a Hamiltonian cycle for all values of $k$. Among various interesting properties of the bipartite Kneser graph $H(n, k)$, we are interested in its automorphism group,  and we want to know how this group acts on its vertex set. Mirafzal [8] determined the automorphism group of $MQ_n = HS(2n,n)= H(2n-1, n-1)$  and  showed that $HS(2n,n)$ is a vertex-transitive non-Cayley graph. Also, he showed that $ H(2n-1, n-1)$ is arc-transitive.   We  show that $H(n, k)$ is a  vertex-and arc-transitive graph for all values of $k$, and hence its connectivity is maximum (Watkins [15]).   Therefore it is a suitable  candidate for some designing in interconnection networks. We   determine the automorphism group  of the graph $H(n, 1)$ and show that $Aut(H(n, 1)) \cong Sym([n] ) \times \mathbb{Z}_2$,  where $\mathbb{Z}_2$ is the cyclic group of order $2$. Also,  we show that  the graph $H(n, 1)$ is a Cayley graph. 

\section{Preliminaries}
In the first step of  our work, we fix some   definitions and notation that we use in
this paper. For all the terminology and notation not defined here, we follow [1,3,4].
In this paper, a graph  $\Gamma = (V, E)$ is a simple, connected and finite graph with vertex set $V$
and edge set $E$. If $v_i$ is adjacent to $v_j$, that is, $\{v_i, v_j \} \in E$, then  we write $v_i \sim v_j$.
Let $\Gamma$ be a graph with automorphism group $Aut(\Gamma)$. The graph $\Gamma$ is called a $vertex-transitive$  graph if  the group $Aut(\Gamma)$ acts transitively  on the vertex set $V$, namely, for any $x, y \in V(\Gamma)$, there is some $\pi$ in $Aut(\Gamma)$,  such that $\pi(x) = y$. For $v\in V(\Gamma)$ and $G=Aut(\Gamma)$, the stabilizer subgroup
$G_v$ is the subgroup of $G$ consisting of all automorphisms that
fix $v$. In the vertex-transitive case all stabilizer subgroups
$G_v $ are conjugate in $G$, and consequently are  isomorphic. The index of $G_v$ in $G$ is given by the equation,  $|G
: G_v| =\frac{|G|}{|G_v|} =|V(\Gamma)|
$.  \  \
 We say that $\Gamma$ is $edge-transitive$ if the group $Aut(\Gamma)$ acts transitively  on the edge set $E$, namely, for any $\{x, y\} ,   \{v, w\} \in E(\Gamma)$, there is some $\pi$ in $Aut(\Gamma)$,  such that $\pi(\{x, y\}) = \{v, w\}$.  We say that $\Gamma$ is $symmetric$ (or $arc-transitive$) if  for all vertices $u, v, x, y$ of $\Gamma$ such that $u$ and $v$ are adjacent,  and also, $x$ and $y$ are adjacent, there is an automorphism $\pi$ in $Aut(\Gamma)$ such that $\pi(u)=x$ and $\pi(v)=y$. We say that $\Gamma$ is $distance-transitive$ if  for all vertices $u, v, x, y$ of $\Gamma$ such that $d(u, v)=d(x, y)$, where $d(u, v)$ denotes the distance between the vertices $u$ and $v$  in $\Gamma$,  there is an automorphism $\pi$ in $Aut(\Gamma)$ such that  $\pi(u)=x$ and $\pi(v)=y.$  It is clear that we have a hierarchy of conditions (Biggs [1]): \

$\textbf{distance-transitive}\Rightarrow \textbf {symmetric}\Rightarrow \textbf {vertex-transitive}$.  \

 The family of distance-transitive graphs includes many interesting and important graphs such as  Hamming graphs, Kneser  graphs and Johnson graphs [1,4].

Let $G$ be a finite group and $\Omega$ a subset of $G$  such that it is closed under taking inverses and does not contain the identity. A $Cayley \ graph$  $\Gamma=Cay(G, \Omega)$ is  a graph whose vertex set and edge set are   as follows:
$V (\Gamma) =G; \,\, \,\ E(\Gamma) = \{\{x, y\} \, | \,\, x^{-1}y \in \Omega\} $
It is well known that every Cayley graph is a vertex-transitive graph [1,4].

 We   show that
$H(n,1)\cong Cay(\mathbb{D}_{2n}, \Omega)$, where  $\mathbb{D}_{2n}=<a,b\  | \  a^n=b^2=1, \   ba=a^{-1}b>$ is the dihedral group of order $2n$ and $\Omega=\{ab, a^{2}b, ... , a^{n-1}b\}$,  which is an  inverse-closed subset of $\mathbb{D}_{2n}-\{1\}$ (note that   $(a^{i}b)^2 =1$ and hence the inverse of $a^{i}b$ is $a^{i}b$).
\section{Main results}
Let $[n]=\{1, 2, ... , n\}$,  and $A$ and $B$ are $m$-subsets of $[n]$. Let $ | A \cap B |=t$ and $\theta$ be  a permutation in $Sym([n])$. It is an easy task to show that  
  $  | f_\theta(A) \cap  f_\theta(B) |=t$, where $ f_\theta (\{x_1, ..., x_m \}) = \{ \theta (x_1), ..., \theta (x_m) \}$.
 Moreover, if $ A\subset B$,  then $  f_\theta(A) \subset f_\theta(B) $. Therefore if $\theta \in Sym[n]$,  then \newline 
$ f_\theta : V( H(n,k) )\longrightarrow V( H(n,k)) $,
 $ f_\theta (\{x_1, ..., x_t \})$ = \newline
$ \{ \theta (x_1), ..., \theta (x_t)\},  t \in \{k,n-k  \} $ \newline
 is an automorphism of $H(n,k) $,  and the mapping \newline
  $ \psi: Sym ([n]) \longrightarrow Aut ( H(n,k) )$, defined by
this  rule  $ \psi ( \theta ) = f_\theta $ is an injection.

 Also, the mapping $\alpha : V(\Gamma)\rightarrow V(\Gamma) $,  defined  by this  rule $\alpha(v) = v^c$,  where
$v^c$ is  the complement of the subset $v$ in  $[n]$, is an automorphism of the bipartite Kneser graph $H(n,k)$, because
if $ A\subset B$, then $ B^c\subset A^c$, and hence if $\{A,B\}$  is an edge of the graph $ H(n,k) $, then $\{\alpha(A),\alpha(B)\}$ is an edge of the graph $ H(n,k) $.\

\begin{lem}
The graph $H(n, k)$ is a vertex-transitive graph.

\end{lem}
\begin{proof}

Let $[n]=\{1, 2, ... , n\}$, $\Gamma=H(n, k)$ and $V = V (\Gamma)$. It is easy to prove that
the graph
$\Gamma$ is a regular bipartite graph. In fact, if
 $V_1=\{ v\in V   | \,\, |v|=k \}$ and $V_2=\{ v\in V   | \,\, |v|=n-k \}$ then $ V=V_1\cup V_2$ and $|V_1| = |V_2 |$=$ {n} \choose{k} $,
 and every edge of $\Gamma$ has a vertex in $V_1$ and  a vertex in $V_2$. Suppose $u, v\in V$.  In the following steps,  we show that $\Gamma$ is  a vertex-transitive graph.

(i)
If both of the two   vertices $u$ and $v$ lie in $V_1$ and $|u \cap v| = t$, then we may assume that
$u = \{ x_1,  ... ,  x_t,  u_1, ... , u_{k-t}\}$ and \newline
$v = \{ x_1,  ... ,  x_t,  v_1, ... , v_{k-t}\}$, where $x_i, u_i, v_i \in [n]$. Let $\sigma$ be a permutation of Sym([n]) such that $\sigma(x_i)=x_i$, $\sigma(u_i)=v_i$ and $\sigma(w_i)=w_i$, where $w_i\in [n]-(u\cup v)$.
Thus, as mentioned early,  $\sigma$ induces an automorphism  $f_{\sigma}:V(\Gamma)\rightarrow V(\Gamma)$
 and we have
$f_{\sigma}\{x_1,  ... ,  x_t,  u_1, ... , u_{k-t}\}=\{ \sigma (x_1),  ... ,  \sigma (x_t), \sigma( u_1), ... , \sigma (u_{k-t})\}$. \newline
 Therefore, 
$f_{\sigma}(u)=v$.

(ii) We now assume that both of the two   vertices $u$ and $v$ lie in $V_2$.
We saw that  the mapping
$\alpha : V(\Gamma)\rightarrow V(\Gamma) $,  defined by the rule  $\alpha(v) = v^c$,  where
$v^c$ is   the complement of the set $v$ in  $[n]$
(for every $v$ in $V$), is  an
automorphism of $\Gamma$. Therefore,   $\alpha(u), \alpha(v)\in V_1$,
 and hence  there is
an automorphism $f_{\sigma}$ in $Aut(\Gamma)$ such that $f_{\sigma}(\alpha(u))=\alpha(v)$, thus
$(\alpha^{-1}f_{\sigma}\alpha)(u)=v$.

 (iii)
Now, let $u\in V_1$ and $v\in V_2$, thus $\alpha(v)\in V_1$,  and hence  there is an automorphism $f_{\sigma}$
in $Aut(\Gamma)$ such that $f_{\sigma}(u)=\alpha(v)$, thus $(\alpha^{-1}f_{\sigma})(u)=v$.
\end{proof}

\begin{thm}
The graph $H(n, k)$ is a symmetric  (or arc-\newline transitive) graph.
\end{thm}
\begin{proof}
Let $\Gamma= H(n, k)$ and $G = Aut(\Gamma)$. Since $\Gamma$ is a vertex-transitive graph, thus it is enough to show that $G_v$ acts transitively on
$N(v)$ for any $v \in V(\Gamma)$ [1, chapter 15]. Without loss of generality, we may assume that  $v=\{ x_1, x_2, ... , x_k\} \in V_1$,   $x_i\in [n]$,  where $V_1$   is as defined in Lemma 3.1. and $k< n-k$. Hence,
 $N(v)=\{ x_1, x_2, ... , x_k, y_1, y_2, ... , y_{n-2k} \} $, where $y_i\in [n]- v$. Suppose $u, w\in N(v)$ such that $|u\cap w|=k+t$,
then we may assume
$u = \{ x_1,  ... ,  x_k,  y_1, ... , y_t, u_1, ... , u_{n-2k-t} \}$ and \newline
$w = \{ x_1,  ... ,  x_k, y_1, ... , y_t, w_1, ... , w_{n-2k-t} \}$. Let $\sigma$ be a permutation of  $Sym([n])$ such that $\sigma(x_i)=x_i$, $\sigma(y_i)=y_i$,  $\sigma(u_i)=w_i$  and $\sigma(z)=z$, where $z\in [n]-(u\cup v \cup w)$.
Therefore,  $\sigma$ induces an automorphism  $f_{\sigma}:V(\Gamma)\rightarrow V(\Gamma)$ and hence we have;\newline
$f_{\sigma}(\{x_1,  ... ,  x_k,  y_1, ... , y_t, u_1, ... , u_{n-2k-t}\})=$ \newline
$ \{ \sigma(x_1),  ... ,  \sigma(x_k),  \sigma(y_1), ... ,\sigma( y_t), \sigma(u_1), ... , \sigma(u_{n-2k-t}) \}$= \newline
$ \{ x_1,  ... ,  x_k, y_1, ... , y_t, w_1, ... , w_{n-2k-t} \} $. 
 Therefore,  $f_{\sigma}(u)=w$.

\end{proof}

 \begin{cor}The connectivity of the  bipartite Kneser graph $H(n, k)$ is
   maximum, namely,  ${n-k}\choose{k}$.

\end{cor}

\begin{proof}\ Since the bipartite Kneser graph $H(n, k)$  is a symmetric
graph,  then it is edge-transitive. On the other hand,  this graph is
 regular  with valency ${n-k}\choose{k}$. Now, note that the connectivity of a
connected edge-transitive graph is equal to its minimum degree (Watkins [15]).

\end{proof}

\

\begin{prop}
 $H(n, 1)$  is a  distance-transitive graph.

\end{prop}

\begin{proof}
Let $[n]=\{1, 2, ... , n\}$ and $v_i= [n]-\{ i \}$, $ i \in [n]$.  We know that  the vertex set of  the bipartite Kneser graph $ \Gamma= H(n,1)$ is $V= \{ \{1\}, \{2\},...,\{n\}, v_1, v_2, ..., v_n   \}$.

Let  $G=Aut(\Gamma)$. We consider the vertex  $v=v_1= \{2,3,...,n  \}$, and  let $ \Gamma_i(v)$ be the set of vertices $w$ in $V(\Gamma)$ at distance $i$ from $v$.  Then $\Gamma_0(v)=\{v\}$,  $\Gamma_1(v)= \{  \{2\},\{3\},...,\{n\} \}$,
$\Gamma_2(v)=\{v_2,v_3,...,v_n \}$ and $\Gamma_3(v)=\lbrace \{1\}\rbrace$. Hence, the diameter of  $\Gamma$  is $3$, because  $\Gamma$ is a vertex-transitive graph. Therefore,
 it is sufficient to show that the vertex-stabilizer $G_v$ is transitive on the set $\Gamma_r(v)$ for every $r\in\{0, 1, 2, 3\}$ [1, chapter 20],  because  $\Gamma$ is a vertex-transitive graph.
  Consider the cycle  $\rho=(2,3,...,n) \in Sym([n])$ and let  $H =<f_{\rho}>$ be the cyclic group generated by $f_{\rho}$ in $G=Aut(\Gamma)$. Then, it is obvious that $ f_{\rho} $  fixes $v$ and the subgroup $H$ is transitive on $\Gamma_1(v)= \{  \{2\},\{3\},...,\{n\} \}$.
 Thus  $G_v$, the stabilizer subgroup of $v$ in $G$, is transitive on the set $\Gamma_1(v)$.  Also, for every $v_i , v_j$ (where \ $ i,j>1$) in $\Gamma_2(v)$,  the transposition $\tau=(i,j)$ in $Sym([n])$ is such that $f_{\tau}$ fixes the vertex $v$ and we have
$ f_{\tau(v_i)}=v_j$. Therefore,    $G_v$ is transitive on the set $\Gamma_2(v)$. It is obvious that the identity subgroup of $Sym([n])$ is transitive on $\Gamma_3(v)$.

\end{proof}

We now want to investigate Cayley properties of the bipartite Kneser graph $H(n,k)$. We show that if $k=1$, then  $H(n,k)$ is a Cayley graph.

\begin{prop}
The bipartite Kneser graph $H(n, 1)$ is a  Cayley graph.

\end{prop}

\begin{proof}
Let $[n]=\{1, 2, ... , n\}$,    $\Delta=H(n, 1)$  and $\Gamma=Cay(\mathbb{D}_{2n}, \Omega)$, where  $\mathbb{D}_{2n}$=$<a,b \  |\, a^n=b^2=1, \  ba=a^{-1}b>$ is the dihedral group of order $2n$,   $\Omega=\{ab, a^{2}b, ... , a^{n-1}b\}$, which is an inverse-closed subset of $\mathbb{D}_{2n}-\{1\}$ (note that $ {(a^ib)}^2=1 $).
 We show that $\Gamma$   is isomorphic to the graph  $H(n, 1)$.
Consider the following  mapping $f$
\begin{equation*}
f: V(\Delta)\longrightarrow V(\Gamma)
\end{equation*}
\begin{equation*}
f(v) =\left\{
\begin{array}{lr}
a^i \,\,\,\,\,\,\,\,\,\,\,\,\,\,\,\ v=\{ i \}, \  i\in [n],   \\

a^{j}b \,\,\,\,\,\,\,\,\,\,\,\,\ v= [n]-\{ j \}, \  j\in [n],  \\
\end{array} \right.
\end{equation*}
It is is clear that $f$ is a  bijective mapping. Let $ \{ i\}$ and $[n]-\{ j  \}$ be two  vertices of $\Delta$,  then 
$$ \{ i \}\sim [n]-\{ j \}\Leftrightarrow  \{ i \} \subset [n]-\{ j \} \,\,\ \Leftrightarrow \,\,\ i \neq j\,\,\ \{ i \}\sim [n]- \{ j \}$$
$$\Leftrightarrow  (a^i)^{-1}a^jb\in \Omega \Leftrightarrow a^i \sim a^jb \Leftrightarrow f(\{ i \})\sim f( [n]- \{ j \}). $$

Note that if $i=j$,  then $(a^i)^{-1}a^jb=b\notin \Omega$. Thus,  $H(n, 1)\cong\Gamma$.
\end{proof} \

  Mirafzal [8] showed that if  $ k \geq 3$, then the bipartite Kneser graph $H(2k-1,k-1)$ is not a Cayley graph. We now ask the following question. \

\

{ \bf Question 1 } For what values of $k$, the bipartite Kneser graph $  H(n,k) $ is a Cayley graph? \

   \

We now proceed  to determine the automorphism group of the graph $H(n,1)$. We know that
 an automorphism of a graph $X$ is a permutation on
its vertex set $V$ that preserves adjacency relations. The automorphism group of $X$, denoted
by $Aut(X)$, is the set of all automorphisms of $X$ with the binary operation of composition of functions. In most situations, it is difficult  to determine the automorphism group
of a graph, but there are various of these results in the literature  and some of the recent works
come in the references  [7,8,9,10,11,13,14].\

\begin{thm}
Let $H(n,1)$ be a bipartite Kneser graph.
 Then, \newline
  $Aut(H(n,1))   \cong Sym([n] ) \times \mathbb{Z}_2$, where $\mathbb{Z}_2$ is the cyclic group of order $2$.

\end{thm}

\begin{proof}
Let $[n]=\{1, 2, ... , n\}$.  Throughout the proof we let $\{i  \}=i$ and $v_i=[n]-\{ i \}$ for  $1 \leq i \leq n$. Hence, the vertex set of  the bipartite Kneser graph $H(n,1)$ is $V= \{ 1, 2,...,n, v_1, v_2, ..., v_n   \}$ and in $H(n,1)$,  $i \sim v_j$, if and only if $ i \neq j$. We now consider the complement of the graph $H(n,1)$. Let $\Gamma$ be the complement of the graph $H(n,1)$. Therefore, the vertex set of $\Gamma$ is $V$ and each edge of $\Gamma$ is of the form $[i,j]$ or $[v_i,v_j]$ for $i,j \in [n],$ when $ i \neq j $ or $[i,v_i]$.
We know that if $X$ is a graph and $ \bar{X}$ is its  complement, then $ Aut(X)= Aut( \bar{X} )  $. Hence, in some cases for determining $Aut(X)$, it is better to find $Aut(\bar{X})$. \

Let  $ G=Aut(\Gamma) $.
 In the first step,  we show that $ | Aut(\Gamma)  | \leq 2(n!)$. \

 Let
$G_1$ be the stabilizer of the vertex $1$ in the graph $\Gamma$.
Let $f\in G_1$. Thus, $f$ is an automorphism of $\Gamma$ that fixes the
vertex 1. Let $N(1)$ be the set of vertices of $\Gamma $ that are adjacent
to the vertex 1. Note that $ N(1)= \{v_1, 2,...,n  \}$. Let $ S= <N(1)> $
be the induced subgraph of $N(1)$ in $\Gamma$. Therefore,  the mapping  $f_{|_{N(1)}}: N(1) \longrightarrow N(1)$, the
restriction of $f$ to $N(1)$,
  is an automorphism of the graph $S$. Since $ v_1= [n]-{1} $ is the unique
isolated vertex in $S$, hence  each automorphism of $S$ fixes the vertex $v_1$.
Since the induced subgraph of $ \{2,...,n  \} $ is a $ (n-1)$-clique in the graph $S$,
thus  if $\theta$ is an element of $Sym(S)$ that fixes $v_1$, then $\theta$
is an automorphism of $S$, and hence $ Aut(S)\cong Sym([n-1]) $. \

We now define the mapping $ \phi : Aut(\Gamma)_1=G_1 \longrightarrow Aut(S)$, by the rule
$\phi (f)= f_{|_{N(1)}}$, for every $f\in Aut(\Gamma)$. It is trivial that $\phi$
is a group homomorphism.  Let $K= Ker(\phi)$, the kernel of the homomorphism $\phi$.
We show that $K= \{ 1 \}$, the identity group. If $g \in K$, then $g(1)=1$ and
 $ g_{|_{N(1)}}=1$, and hence $g(v_1)=v_1$,  $g(j)=j$ for  $(2\leq j \leq n)$. For each $j$ for $1 \leq j \leq n $, vertices
 $v_j$ and $j$ are adjacent, therefore  vertices $g(v_j)$ and $g(j)=j$ are adjacent,  and it follows that
 $g(v_j)=v_j$. Therefore,  we must have $g=1$, the identity automorphism of  $ \Gamma $,  and consequently
 $Ker(\phi)= \{ 1\}$, and thus $\phi$ is an injection. \

 Since $\phi$ is an injection, then we have 
 \newline $ G_1 \cong \phi(G_1) \leq Sym([n-1]),  $  \newline 
 where $G_1= Aut (\Gamma)_1$, and hence
 $| G_1| \leq (n-1)!$.
On the other hand, $\Gamma$ is a vertex-transitive graph, hence by the orbit-stabilizer theorem,
we have $|V(\Gamma)|=2n= \frac{|G|}{|G_1|}$,  which implies that\

\

\centerline{$ |G|=2n|G_1|\leq 2n((n-1))! =2(n!) $ \ \ \ \ \ \ \ \ (1) } \

  Since $  \Gamma$ is the complement of the graph $H(n,1) $,  we now conclude that 
 $|  Aut(H(n,1)) | \leq  2(n!) $. \newline
Now,  if we find  a subgroup $H$ of $G$ of order $2(n!)$,   where $G$ is  the automorphism group
 of $H(n,1)$,  then from (1) we can conclude
  that $G =H$.\

 We know that the set  $ \{ f_\theta |\  \theta \in Sym ([n]) \} = K $,  is a subgroup of $ Aut ( H(n,1)) $.  Also, we know that the  mapping  $ \alpha : V( \Gamma ) \longrightarrow V(\Gamma), \  \alpha( v ) = v^c $  where $ v^c $ is the complement of the set $ v $ in $[n]$, is
 also an automorphism of $ H(n,1)$,  namely,   $\alpha \in G = Aut (H(n,1 )) $. \

We show that $ \alpha \not\in K $. If $ \alpha \in K $,  then there is a $ \theta \in Sym ([n]) $ such that $ f_\theta = \alpha $. Since $ o( \alpha ) = 2 $ $( o(\alpha )= $ order of $ \alpha ) $,  then $ o(f_\theta ) = o( \theta ) = 2 $. We assert that
$ \theta $ has no fixed points,  namely,  $ \theta (x) \neq x $, for every  $ x \in [n] $. In fact,  if $ x \in [n] $, and $ \theta (x) = x $, then for the $ k $-set $ v = \{x, y_1, ..., y_{k-1} \} \subseteq [n] $, we have 
\newline $ f_{\theta }( v ) = \{ \theta (x), \theta (y_1), ..., \theta( y_{k-1}) \} = \{x, \theta(y_1), ..., \theta (y_{k-1}) \} $,  \newline
and  hence  $ x \in f_\theta (v) \cap v $, and therefore  $ f_\theta (v) \neq v^c = \alpha (v) $,   which is a contradiction.
Therefore,  $ \theta $ has a form such as   $ \theta = ( x_1, y_1 )...(x_m, y_m)$,  where $ ( x_i, y_i)$ is a transposition of $Sym([n]) $, which is impossible if $n $ is an odd  integer (note that $2m=n$). We now assume that $n=2m$ is an even integer. Then,  for the $m$-subset $ v= \{x_1,y_1, x_2,...,x_{m-1}  \}$ of $[n]$,  we have
\

\

   \centerline{$ \alpha (v) = f_{ \theta }(v)= \{ \theta (x_1), \theta (y_1), \theta ( x_{m-1} )   \} $=
     $ \{ y_1, x_1,..., \theta (x_{m-1})   \},  $} \

 and thus  $ x_1,y_1 \in f_ { \theta} (v) \cap v $,   hence $ f_\theta (v) \neq v^c = \alpha (v) $, which is also  a contradiction. \

We assert that  for every $\theta \in Sym([n])$, we have $f_{\theta}\alpha= \alpha f_{\theta}$. In fact, if $v=\{ x_1,...,x_k \}$ is a $k$-subset of $[n]$, then there are $y_j \in [n]$ for  $  1\leq j \leq n-k $,  such that $[n]= \{x_1,...,x_k, y_1,...,y_{n-k }   \}$. Now we have 
 $f_{\theta}\alpha(v)=
  f_{\theta} \{ y_1,...,y_{n-k } \}= \{ \theta(y_1), ...,  \theta(y_{n-k })   \}$. \newline  
 On the other hand, we have
\newline $ \alpha f_{\theta}(v) =\alpha \{ \theta(x_1),..., \theta(x_k) \}= \{ \theta(y_1), ...,  \theta(y_{n-k})   \}, $ \newline
because $[n] = \theta([n])= \{\theta(x_1),...,\theta(x_k), \theta(y_1),...,\theta(y_{n-k} )        \}$.  Consequently,   $f_{\theta}\alpha(v)=  \alpha f_{\theta}(v)  $.   We now deduce that $f_{\theta}\alpha=  \alpha f_{\theta}  $. \

Note that if $X$ is a group and $ Y,Z$ are subgroups of $X$, then  the subset $YZ=\{  yz \ |  \  y\in Y , z\in Z \}$ is a subgroup in X if and only if $YZ=ZY$.  According to this fact,    we conclude that  $ K  < \alpha > $ is a subgroup of $G$.

         Since $ \alpha \not\in K $ and $ o( \alpha ) =2$,  then $ K < \alpha > $ is a subgroup of $ G $ of order\

\

\centerline{ $ \frac { |K| |< \alpha > |}{ | K \cap < \alpha > |} = 2|K| = 2( n )! $ } \

Now,  since  by (1) $ | G | \leq 2( n! ) $,  then $ G =  K  < \alpha > $. On the other hand,  since $f_{\theta}\alpha=  \alpha f_{\theta}  $, for every $ \theta \in Sym([n]) $, then  $ K $ and
$< \alpha >$ are normal subgroups of $ G $. Thus,   $ G $    is a direct product of two  groups $K$ and $< \alpha >$,
  namely, we have  $ G = K \times < \alpha > \cong Sym ([n])\times \mathbb{Z}_2 $.\

\end{proof}
 Mirafzal  [8] proved the following theorem.

\begin{thm} $[8]$
Let $n=2k-1$. Then, for the    bipartite Kneser graph $H(n,k-1)$, we have 
  $Aut(H(n,k)) \cong Sym([n]) \rtimes \mathbb{Z}_2$, where $\mathbb{Z}_2$ is the cyclic group of order $2$.
\end{thm}

We can show by Theorem 3.7.  and   the last part of the proof of Theorem 3.6.  the following theorem.\
\begin{thm} 
Let $n=2k-1$. Then, for the    bipartite Kneser graph $H(n,k-1)$, we have 
  $Aut(H(n,k)) \cong Sym([n])  \times \mathbb{Z}_2$, where $\mathbb{Z}_2$ is the cyclic group of order $2$.

\end{thm}

Now, it is natural to ask  the following question. \

\

{ \bf Question 2 } Is the above theorem true for all possible values of $n,k$ ( $2k \neq  n$)? \

\

{\bf Acknowledgements}
\

\

The authors are thankful to the anonymous reviewer for his (her) valuable comments and suggestions.

\end{document}